\documentclass[a4paper,12pt]{amsart}
\usepackage{amsmath, amssymb, amsthm}
\usepackage{verbatim}
\usepackage{hyperref}

\newcounter{theorem}

\newtheorem{prop}[theorem]{Proposition}

\theoremstyle{remark}
\newtheorem*{remark*}{Remark}

\newcommand{\tens}{\otimes}

\newcommand{\Cu}{\mathcal{C}u}

\newcommand{\diag}{\mathrm{diag}}

\newcommand{\labelledthing}[2]{\hspace{4pt}\buildrel {#2} \over #1 \hspace{3pt}} 

\newcommand{\labelledrightarrow}{\labelledthing{\longrightarrow}}

\newcommand{\X}{\mathbb{X}}
\newcommand{\Y}{\mathbb{Y}}

\begin{document}
\title{Corrigendum to ``Regularity for stably projectionless, simple $C^*$-algebras''}
\author{Henning Petzka}
\address{\hskip-\parindent
Henning Petzka, Mathematisches Institut der WWU M\"unster, Einsteinstra\ss{}e 62, 48149 M\"unster, Germany.}
\email{petzka@uni-muenster.de}
\author{Aaron Tikuisis}
\address{\hskip-\parindent Aaron Tikuisis, Institute of Mathematics, School of Natural and Computing Sciences, University of Aberdeen, AB24 3UE, Scotland.}
\email{a.tikuisis@abdn.ac.uk}
\thanks{Research partially supported by EPSRC (grant no.\ EP/N002377/1), NSERC (PDF, held by AT), and by the DFG (SFB 878).}


\newcommand{\cstar}{$C^*$}

\begin{abstract}
An error is identified and corrected in the construction of a non-$\mathcal Z$-stable, stably projectionless, simple, nuclear \cstar-algebra carried out in a paper by the second author.
\end{abstract}

\maketitle

\section*{The problem}

The construction in Section 4 of the second author's paper \cite{ProjlessReg}, used to prove \cite[Theorem 4.1]{ProjlessReg}, contains a vital error.
The construction is meant to produce a simple \cstar-algebra with perforation in its Cuntz semigroup, as an inductive limit of stably projectionless subhomogeneous \cstar-algebras.

The notation set out in \cite{ProjlessReg} will be reused here, mostly without recalling the definitions.

The idea is to use generalized Razak building blocks $R(\mathbb X,k) \subseteq C(X,M_{k+1})$ (as defined in \cite[Section 4.2]{ProjlessReg}) as the stably projectionless building blocks of the inductive system; the connecting maps are unitary conjugates of restrictions of diagonal maps $D_{\alpha_1,\dots,\alpha_p}:C(X,M_n) \to C(Y,M_m)$ (as defined in \cite[Section 4.1]{ProjlessReg}).

For generalized Razak building blocks $R(\mathbb X,k) \subseteq C(X,M_{k+1})$ and $R(\mathbb Y,\ell) \subseteq C(X,M_{\ell+1})$, \cite[Proposition 4.3]{ProjlessReg} characterizes when a diagonal map $D_{\alpha_1,\dots,\alpha_p}:C(X,M_{k+1}) \to C(Y,M_{\ell+1}) \otimes M_m$ is unitarily conjugate to a map which sends $R(\mathbb X,k)$ into $R(\mathbb Y,\ell) \otimes M_m$.
The characterization includes the equations
\begin{align}
\label{eq1}
ka_0 + (k+1)a_1 &= (m-s(k+1))\ell, \text{ and} \\
\label{eq2}
kb_0 + (k+1)b_1 &= (m-s(k+1))(\ell+1),
\end{align}
where $a_0,a_1,b_0,b_1$, and $s$ count certain values of the maps $\alpha_1,\dots,\alpha_p$; they additionally satisfy
\begin{equation}
\label{eq3}
p = a_0+a_1+s\ell = b_0+b_1+s(\ell+1).
\end{equation}
In \cite[Remark 4.4]{ProjlessReg}, a specific (parametrized) solution is provided to the condition in \cite[Proposition 4.3]{ProjlessReg}, and this solution is used in \cite[Section 4.4]{ProjlessReg} to construct the example.

Implicit in the definition of diagonal maps in \cite[Section 4.1]{ProjlessReg} is that they are unital (as maps $C(X,M_n) \to C(Y,M_m)$).
In the case of \cite[Proposition 4.3]{ProjlessReg}, this means that 
\begin{equation}
\label{eq4}
p(k+1)=m(\ell+1).
\end{equation}
However, the solution provided in \cite[Remark 4.4]{ProjlessReg} does not satisfy \eqref{eq4}.
In fact, some algebraic manipulation of the equations in \cite[Proposition 4.3]{ProjlessReg} shows that there are not very many solutions at all.
Certainly, suppose that $m,\ell,p,s,a_0,a_1,b_0,b_1$ satisfy \eqref{eq1}, \eqref{eq2}, \eqref{eq3}, and \eqref{eq4}.
Combining \eqref{eq3} and \eqref{eq4} yields 
\[ (b_0+b_1+s(\ell+1))(k+1) = m(\ell+1). \]
Subtracting \eqref{eq2} from this produces $b_0=0$.
Likewise, one obtains $a_0=m$.

Crucial to the construction in \cite{ProjlessReg} is the use of both coordinate projections and flipped coordinate projections among the eigenmaps in the diagonal map $D_{\alpha_1,\dots,\alpha_p}$.
As intimated in \cite[Remark 4.4]{ProjlessReg}, there may be up to $\max\{a_0,b_1\}$ coordinate projections and $\max\{a_1,b_0\}$ flipped coordinate projections.
To get perforation, the number of coordinate projections and flipped coordinate projections needs to be a very large fraction of the total number of eigenmaps.
Since solutions to \cite[Proposition 4.3]{ProjlessReg} necessarily have $b_0=0$, it is actually not possible to get perforation in the Cuntz semigroup with this kind of construction.

\section*{The solution}

Here we describe a correction to the construction in \cite[Section 4]{ProjlessReg}, permitting a correct proof of \cite[Theorem 4.1]{ProjlessReg}.
The solution is to allow slightly more general diagonal maps which include some copies of the zero representation.

Let $X,Y$ be compact Hausdorff spaces and let $\alpha_1,\dots,\alpha_p:Y \to X$ be continuous functions.
Suppose that $m,n,r \in \mathbb N$ satisfy $np+r=m$.
Define $D_{\alpha_1,\dots,\alpha_p;r}:C(X,M_n) \to C(Y,M_m)$ by
\begin{align*}
D_{\alpha_1,\dots,\alpha_p;r}(f) & := \diag(f \circ \alpha_1, f \circ \alpha_2, \dots, f \circ \alpha_p, 0_r) \\
&:= \left( \begin{array}{cccc}
f \circ \alpha_1 & 0 & \cdots & 0 \\
0 & f \circ \alpha_2 & \ddots & \vdots \\
\vdots & \ddots & \ddots & 0 \\
0 & \cdots & f \circ \alpha_p & 0 \\
0 & \cdots & 0 & 0_r
\end{array}\right),
\end{align*}

We have the following generalization of \cite[Proposition 4.2]{ProjlessReg} (the only difference being that the map $D_{\alpha_1^{(i)},\dots,\alpha_{p_i}^{(i)}}$ is replaced by the more general $D_{\alpha_1^{(i)},\dots,\alpha_{p_i}^{(i)};r_i}$).
The proof is exactly the same.

\begin{prop}\label{DiagonalPerforation}
Let
\[ A_1 \labelledrightarrow{\phi_1^2} A_2 \labelledrightarrow{\phi_2^3} \cdots \]
be an inductive limit, such that for each $i$, the algebra $A_i$ is a subalgebra of $C(X_i,M_{m_i})$ and $\phi_i^{i+1} = Ad(u) \circ D_{\alpha_1^{(i)},\dots,\alpha_{p_i}^{(i)};r_i}$ for some unitary $u \in C(X_{i+1},M_{m_{i+1}})$ (so that $m_{i+1}=m_ip_i+r_i$).
Suppose that $X_i$ contains a copy $Y_i$ of $[0,1]^{d_1\cdots d_{i-1}}$ such that
\begin{itemize}
\item $A_i|_{Y_i} = C(Y_i,M_{m_i})$,
\item for $t=1,\dots,d_i$, $\alpha_t^{(i)}|_{Y_{i+1}}$ takes $Y_{i+1}$ to $Y_i$ via the $t^\text{th}$ coordinate projection $([0,1]^{d_1 \cdots d_{i-1}})^{d_i} \to [0,1]^{d_1 \cdots d_{i-1}}$, and
\item for $t=d_i+1,\dots,p_i$, $\alpha_t^{(i)}|_{Y_{i+1}}:Y_{i+1} \to X_i$ factors through the interval.
\end{itemize}
If
\[ \prod_{i=1}^\infty \frac{d_{i+1}}{p_i} > 0 \]
and $p_i > 1$ for all $i$ then for any $n \in \mathbb{N}$, there exists $[a],[b] \in \Cu(\varinjlim A_i)$ and $k \in \mathbb{N}$ such that
\[ (k+1)[a] \leq k[b] \]
yet $[a] \not\leq n[b]$.
\end{prop}

We have the following generalization of \cite[Proposition 4.3]{ProjlessReg}; the diagonal map $D_{\alpha_1,\dots,\alpha_p}$ of \cite[Proposition 4.3]{ProjlessReg} is replaced by the more general $D_{\alpha_1,\dots,\alpha_p;r}$.
This results in a looser condition in (ii) (compare \eqref{eq1}, \eqref{eq2} to \eqref{eq1b}, \eqref{eq2b} respectively). 
The proof is nearly the same and contains no new tricks.

\begin{prop}\label{CharacterizingExistence-Razak}
Let $\X = (X,x_0,x_1), \Y = (Y,y_0,y_1)$ be double-pointed spaces and let $k,\ell,m,p,r$ be natural numbers such that
\begin{equation}
\label{eq4b}
p(k+1)+r=m(\ell+1).
\end{equation}
Let $\alpha_1,\dots,\alpha_p:Y \to X$ be continuous maps.
Then the following are equivalent:
\begin{enumerate}
\item There exists a unitary $u \in C(Y,M_{\ell+1}) \tens M_m$ such that 
\[ uD_{\alpha_1,\dots,\alpha_p;r}(R(\X,k))u^* \subseteq R(\Y,\ell) \tens M_m; \text{ and} \]
\item Counting multiplicity we have
\begin{align*}
\{\alpha_1(y_0),\dots,\alpha_p(y_0)\} &= a_0\{x_0\} \cup a_1\{x_1\} \cup \ell\{z_1\} \cup \cdots \cup \ell\{z_s\} \text{ and} \\
\{\alpha_1(y_1),\dots,\alpha_p(y_1)\} &= b_0\{x_0\} \cup b_1\{x_1\} \cup (\ell+1)\{z_1\} \cup \cdots \cup (\ell+1)\{z_s\}
\end{align*}
for some points $z_1,\dots,z_s \in X$, and some natural numbers $a_0,a_1,b_0,b_1$ satisfying
\begin{align}
\label{eq1b}
ka_0 + (k+1)a_1 &= (m-s(k+1)-q)\ell, \text{ and} \\
\label{eq2b}
kb_0 + (k+1)b_1 &= (m-s(k+1)-q)(\ell+1),
\end{align}
for some $q \in \mathbb N$.
\end{enumerate}
\end{prop}

Here is a solution to \eqref{eq3}, \eqref{eq4b}, \eqref{eq1b}, and \eqref{eq2b}, parametrized by $s,k,u \in \mathbb N_{>0}$; it is almost the same as the solution in \cite[Remark 4.4]{ProjlessReg} with the notable difference of being correct.
\begin{align*}
\ell &:= k+1+2u, \\
m &:= (k^2+3k+1)s, \\
a_0 &:= (k+1)(k+1+u)s, \quad a_1 := ksu, \\
b_0 &:= (k+1)su, \quad b_1 := k(k+2 + u)s, \\
r &:= (k^2+2k+ku-u)s, \\
q &:= ks, \\
p &:= (k^2 + 2ku + 3k + 3u + 2)s.
\end{align*}

The construction in \cite[Section 4.4]{ProjlessReg} proceeds using this solution in place of the one in \cite[Remark 4.4]{ProjlessReg}.
In essence, the only difference is that the assignment
\[ m_{i+1} := m_i(k_i+1)^2s_i \]
is replaced by
\[ m_{i+1} := m_i(k_i^2+3k_i+1)s_i. \]

As opposed to the original (though incorrect) construction in \cite{ProjlessReg}, it is not obvious that the algebra $A$ constructed with these corrections has a tracial state (as opposed to only having a densely defined trace).
One need not be concerned that this causes problems in proving the desired properties of this example, since nowhere in the statement or proof of \cite[Theorem 4.1]{ProjlessReg} (nor elsewhere in \cite{ProjlessReg}) is it used that $A$ has a tracial state.

This correction thereby provides a proof of \cite[Theroem 4.1]{ProjlessReg}.

\end{document}